\documentclass[12pt]{article}
\usepackage[dvips]{graphicx}
\usepackage{amsmath,amssymb,amsthm,bm}
\bmdefine{\BX}{{\bm X}}
\bmdefine{\Bx}{{\bm x}}
\bmdefine{\By}{{\bm y}}
\bmdefine{\BY}{{\bm Y}}
\bmdefine{\Bz}{{\bm z}}
\bmdefine{\Ba}{{\bm a}}
\bmdefine{\Bd}{{\bm d}}
\bmdefine{\Bu}{{\bm u}}
\bmdefine{\Bm}{{\bm m}}
\bmdefine{\Bi}{{\bm i}}
\bmdefine{\Bc}{{\bm c}}
\bmdefine{\BT}{{\bm T}}
\bmdefine{\Bp}{{\bm p}}
\bmdefine{\Bs}{{\bm s}}
\bmdefine{\Btheta}{{\bm \theta}}
\bmdefine{\Bmu}{{\bm \mu}}
\bmdefine{\Bt}{{\bm t}}
\bmdefine{\Blambda}{{\bm \lambda}}
\bmdefine{\Bbeta}{{\bm \beta}}
\bmdefine{\Bpsi}{{\bm \psi}}
\bmdefine{\Bpi}{{\bm \pi}}
\bmdefine{\Bv}{{\bm v}}
\bmdefine{\Bzero}{{\bm 0}}
\bmdefine{\Bone}{{\bm 1}}
\bmdefine{\Bvarepsilon}{{\bm \varepsilon}}

\bmdefine{\BI}{{\rm {\bm I}}}
\bmdefine{\BV}{{\rm {\bm V}}}

\newtheorem{theorem}{Theorem}[section]
\newtheorem{proposition}[theorem]{Proposition}

\newtheorem{corollary}[theorem]{Corollary}

\theoremstyle{definition}
\newtheorem{example}[theorem]{Example}
\newtheorem{definition}[theorem]{Definition}

\title{
Characterizations of indicator functions and contrast
representations of 
 fractional factorial designs
 with multi-level factors
}
\author{Satoshi Aoki%
\thanks{Graduate School of Science, Faculty of Science, Kobe University.}
}
\date{}

\begin{document}
\maketitle

\begin{abstract}
\noindent
A polynomial indicator function of designs is first 
introduced by Fontana, Pistone and Rogantin (2000) for two-level
designs. They give the structure of the indicator function of
 two-level designs, especially from the viewpoints of the orthogonality
 of the designs. Based on the structure, they use the indicator
 functions to classify all the orthogonal fractional factorial designs
 with given sizes using computational algebraic software.
In this paper, generalizing the results on two-level
 designs, the structure of the indicator functions for multi-level
 designs is derived. We give a system of algebraic equations for
 the coefficients of indicator functions of fractional factorial designs
 with given orthogonality. We also give another representation of the
 indicator function, a contrast representation, which reflects the size
 and the orthogonality of the corresponding design directly. 
The contrast representation is
determined by a contrast matrix, and does
 not depend on the level-coding, which is one of the advantages of it. 
We use
 these results to classify orthogonal $2^3\times 3$ designs with
 strength $2$ and orthogonal  $2^4\times 3$ designs with strength $3$ by
a computational  algebraic software.\\

\end{abstract}

\section{Introduction}
\label{sec:intro}
Applications of Gr\"obner basis theory to various problems of statistics
arises in early 1990s. One of the first works in this
developing field, a computational algebraic statistics, is given by
Pistone and Wynn (\cite{Pistone-Wynn-1996}), where the Gr\"obner
basis theory is applied to 
the identifiability problem in the design of experiments. 
After this work, various algebraic
approaches to the problems in the design of experiments 
are presented by researchers both in the fields of algebra and
statistics. A theory
of the indicator function of fractional factorial designs 
is one of the early results in this branch.

The indicator function is first introduced by Fontana, Pistone and
Rogantin (\cite{Fontana-Pistone-Rogantin-2000}) 
for two-level fractional
factorial designs. In \cite{Fontana-Pistone-Rogantin-2000}, based on the
results of 
\cite{Pistone-Wynn-1996}, 
one-to-one correspondence between the design and its 
indicator function is shown.
This correspondence enables us to translate various statistical
concepts to algebraic concepts, i.e., 
various results on the fractional factorial designs
can be interpreted to the structure of their indicator functions. 
For example, abberation and resolution are important concepts in design
of experiments, and 
there is a well-established history starting with
\cite{Box-Hunter-1961} for two-level fractional factorial designs. 
An important contribution of \cite{Fontana-Pistone-Rogantin-2000} is to
characterize these concepts as the structure
of the indicator functions. 

To illustrate the motivation of this paper, we 
glance at the arguments of \cite{Fontana-Pistone-Rogantin-2000} by
examples. Note that the necessary definitions on the designs and indicator
functions will be given in Section
\ref{sec:definition-indicator-function}. 
Figure \ref{fig:example-D1-D2} shows 
examples of two-level fractional factorial designs.
\begin{figure*}[htbp]
\[
\begin{array}{|rrrrr|}
\multicolumn{5}{l}{F_1}\\
\multicolumn{1}{c}{x_1} & x_2 & x_3 & x_4 & \multicolumn{1}{c}{x_5}\\ \hline
 1 & 1 & 1 & 1 & 1\\
 1 & 1 &-1 & 1 &-1\\
 1 &-1 & 1 &-1 & 1\\
 1 &-1 &-1 &-1 &-1\\
-1 & 1 & 1 &-1 &-1\\
-1 & 1 &-1 &-1 & 1\\
-1 &-1 & 1 & 1 &-1\\
-1 &-1 &-1 & 1 & 1\\ \hline
 \end{array}
\hspace*{10mm}
\begin{array}{|rrrr|}
\multicolumn{4}{l}{F_2}\\
\multicolumn{1}{c}{x_1} & x_2 & x_3 & \multicolumn{1}{c}{x_4}\\ \hline
 1 & 1 & 1 & 1\\
 1 & 1 &-1 &-1\\
 1 &-1 & 1 &-1\\
-1 & 1 &-1 &-1\\
-1 &-1 & 1 &-1\\
-1 &-1 &-1 & 1\\ \hline
\multicolumn{4}{c}{}\\
\multicolumn{4}{c}{}
 \end{array}
\]
\caption{Examples of fractional factorial designs for two-level
 factors. Left($F_1$): a regular fractional factorial
 design with the defining
 relation $x_1x_2x_4 = x_1x_3x_5 = 1$. Right($F_2$): an example of
 nonregular designs.}
\label{fig:example-D1-D2}
\end{figure*}
We code the levels of each factor as $\{-1,1\}$ according to 
\cite{Fontana-Pistone-Rogantin-2000}. 
For each design, each row of the table shows the combination of the
levels of the factors $x_i$'s for each experimental run, and each column
corresponds to each
factor. For example, the design $F_1$ is a fractional factorial
design for $5$ two-level factors $x_1,\ldots,x_5$, composed of $8$
points in $\{-1,1\}^5$,
\[
 \{(1,1,1,1,1),(1,1,-1,1,-1),\ldots,(-1,-1,-1,1,1)\}.
\]
In the field of design of experiments, 
$F_1$ is known as a regular fractional factorial design with the
defining relation $x_1x_2x_4 = x_1x_3x_5 = 1$. On the other hand, the design
$F_2$ is an example of nonregular designs. For details on
the regularity of designs, see \cite{Wu-Hamada-2009} for example. 

The indicator functions of $F_1$ and $F_2$ are given as follows, respectively.
\[\begin{array}{cl}
F_1: & \displaystyle 
f_1(x_1,x_2,x_3,x_4,x_5) = \frac{1}{4} + \frac{1}{4}(x_1x_2x_4 +
 x_1x_3x_5 +
  x_2x_3x_4x_5)\\
F_2: & \displaystyle f_2(x_1,x_2,x_3,x_4) = \frac{3}{8} - \frac{1}{8}x_4 +
 \frac{1}{8}(x_1x_2 + x_1x_3 - x_2x_3)
+ \frac{1}{8}(x_1x_3x_4 +
 x_2x_3x_4) + \frac{3}{8}x_2x_3x_4
\end{array}
 \]
We see, for example, $f_1(x_1,\ldots,x_5) = 1$ for $8$ points in $F_1$,
and $f_1(x_1,\ldots,x_5) = 0$ for the other $24$ points not in $F_1$.  
The indicator function of the design of $n$ two-level 
factors, $x_1, x_2,\ldots,x_n$ has a unique polynomial representation 
of the form
\begin{equation}
 f(x_1,\ldots,x_n) = \sum_{\Ba \in \{0,1\}^n}\theta_{\Ba}\Bx^{\Ba}, 
\label{eqn:indicator-function-two-level}
\end{equation}
where $\Bx^{\Ba} = \displaystyle\prod_{i = 1}^nx_i^{a_i}$ and
$\Ba = (a_1,\ldots,a_n) \in \{0,1\}^n$.
As is shown in \cite{Fontana-Pistone-Rogantin-2000}, 
the set of the coefficients $\{\theta_{\Ba}\}_{\Ba \in \{0,1\}^n}$ has 
all the information of the corresponding design. 
For example, we see the following facts from the coefficients of the
indicator functions $f_1$ and $f_2$ for $F_1$ and $F_2$.
\begin{itemize}
\item The constant term $\theta_{(0,\ldots,0)}$ shows the ratio 
between the size of the design to the size of the full factorial
      design. In fact, $F_1$ is a $1/4$ fraction of the full factorial
      $2^5$ design, and $F_2$ is a $3/8$ fraction of the full factorial
      $2^4$ design.
\item The coefficient of the main effect term $\theta_{\Ba}$, $\sum_j
      a_j = 1$, i.e., the coefficient of the monomial $\Bx^{\Ba}$ with
      the degree $1$, shows the
      ``balance'' of two levels for this factor. In fact, for $F_1$, 
      $\theta_{(1,0,0,0,0)} = \cdots = \theta_{(0,0,0,0,1)} = 0$ shows
      $F_1$ is an equireplicated design, i.e., two levels
      appear equally often
      for each factor. On the other hand, for $F_2$, $\theta_{(1,0,0,0)}
      = \theta_{(0,1,0,0)} = \theta_{(0,0,1,0)} = 0$ and
      $\theta_{(0,0,0,1)} \neq 0$ show $F_2$ is
      equireplicated for factors $x_1,x_2,x_3$ but not for $x_4$.
\item The coefficient of the two-factor interaction term
      $\theta_{\Ba}$, $\sum_j a_j = 2$, i.e., the coefficient of the
      monomial $\Bx^{\Ba}$ with the degree $2$, 
      shows the ``orthogonality'' of the design. In fact, for $F_1$, 
      $\theta_{(1,1,0,0,0)} = \cdots = \theta_{(0,0,0,1,1)} = 0$ 
      shows $F_1$ is an orthogonal design, i.e., 
      possible combinations of levels, $(-1,-1), (-1,1), (1,-1), (1,1)$,
      appear equally often for each pair of the factors. On the other
      hand, for $F_2$, $\theta_{(1,0,0,1)} = \theta_{(0,1,0,1)} =
      \theta_{(0,0,1,1)} = 0$ shows that the factor $x_4$ is orthogonal to
      each of the other factors, whereas $\theta_{(1,1,0,0)},
      \theta_{(1,0,1,0)}, \theta_{(0,1,1,0)} \neq 0$ shows that $x_1,
      x_2, x_3$ are not
      orthogonal in each other. 
\end{itemize}
In other words, statistical concepts such as aberration and
resolution can be
related to the structure of the corresponding indicator functions
directly for two-level designs. 
In particular, the structure of the indicator function of
regular two-level designs can be characterized by their defining
relations, and are fully revealed. See \cite{Ye-2003} for
detail. Another characterization of the indicator function of two-level
designs relating the $D$-optimality of the design 
is given by the author in \cite{Aoki-Takemura-2009}. 

In \cite{Fontana-Pistone-Rogantin-2000}, these structures of 
the indicator
function are applied to the classification of the design, which is also
the object of this paper. The argument of
\cite{Fontana-Pistone-Rogantin-2000} is as follows. For the indicator
function (\ref{eqn:indicator-function-two-level}) of two-level designs,
the set of the coefficients $\{\theta_{\Ba}\}_{\Ba \in \{0,1\}^n}$
satisfies a system of algebraic equations
\begin{equation}
 \theta_{\Ba} = \sum_{\Ba' \in \{0,1\}^n}\theta_{\Ba'}\theta_{\Ba +
 \Ba'},\ \ \Ba \in \{0,1\}^n,
\label{eqn:relation-theta-2-level}
\end{equation}
where the sum $\Ba + \Ba'$ is considered under ``mod $2$'' (Proposition 3.7
of \cite{Fontana-Pistone-Rogantin-2000}). Therefore, adding 
constraints for some orthogonality of the designs 
to (\ref{eqn:relation-theta-2-level}), we have a system of
algebraic equations 
having the designs with these orthogonality as the solutions.
For example, for the case of $n = 5$, additional constraints
\begin{equation}
 \theta_{(0,0,0,0,0)} = \frac{1}{4},\ \ 
 \theta_{\Ba} = 0\  \mbox{for}\ \sum_j a_j = 1,2
\label{eqn:example-constraints-2-level}
\end{equation}
to (\ref{eqn:relation-theta-2-level}) yields a system of algebraic
equations having all the orthogonal designs with the size $8$ as the
solution (and $F_1$ corresponds to one of the solutions). 
In this way, the complete lists of the orthogonal designs for
$n = 4,5$ are computed by a computational algebraic software in
\cite{Fontana-Pistone-Rogantin-2000}. Recall that
solving a system of algebraic equations is a fundamental problem where
the theory of Gr\"obner basis is used.

In this paper, we consider generalization of the above argument on 
two-level designs to
general fractional factorial designs. Note that the direct
relations between the size and orthogonality of 
designs and their indicator functions 
are obtained only for two-level designs. 
To see this, consider a fractional factorial design of three-level
factors $F_3$ displayed in Figure \ref{fig:example-F3}.
$F_3$ is 
a regular fractional factorial design with the defining relation $x_1 +
x_2 + x_3 = x_1 +
2x_2 + x_4 = 0\ ({\rm mod}\ 3)$. 
\begin{figure*}[htbp]
\[
\begin{array}{|rrrr|}
\multicolumn{4}{l}{F_3}\\
\multicolumn{1}{c}{x_1} & x_2 & x_3 & \multicolumn{1}{c}{x_4}\\ \hline
-1 &-1 &-1 & 0\\
-1 & 0 & 1 & 1\\
-1 & 1 & 0 &-1\\
 0 &-1 & 1 &-1\\
 0 & 0 & 0 & 0\\
 0 & 1 &-1 & 1\\
 1 &-1 & 0 & 1\\
 1 & 0 &-1 &-1\\
 1 & 1 & 1 & 0\\ \hline
 \end{array}
\]
\caption{A regular fractional factorial design for three-level
 factors with the defining
 relation $x_1 + x_2 + x_3 = x_1 + 2x_2 + x_4 = 0\ ({\rm mod}\ 3)$.}
\label{fig:example-F3}
\end{figure*}
Though $F_3$ is a regular design (and therefore the resolution of $F_3$
is seen in its defining relation), the structure of its indicator
function seems complicated as follows.
\begin{equation}
\begin{array}{cl}
\multicolumn{2}{l}{f(x_1,x_2,x_3,x_4)}\vspace*{2mm}\\
 = & \displaystyle 1 - x_1^2 - x_2^2 - x_3^2 - x_4^2
+ x_1^2x_2^2 + x_1^2x_3^2 + x_1^2x_4^2
+ x_2^2x_3^2 + x_2^2x_4^2 + x_3^2x_4^2\vspace*{2mm}\\
& + \displaystyle \frac{1}{4}(x_1^2x_2x_3 - x_1^2x_2x_4 + x_1^2x_3x_4 +
 x_1x_2^2x_3 +
 x_1x_2^2x_4 - x_2^2x_3x_4\\
&\hspace*{7mm} {} + x_1x_2x_3^2 - x_1x_3^2x_4 + x_2x_3^2x_4 -
 x_1x_2x_4^2 - x_1x_3x_4^2 - x_2x_3x_4^2)\vspace*{2mm}\\
& - \displaystyle \frac{3}{4}(x_1^2x_2^2x_3^2 + x_1^2x_2^2x_4^2 +
 x_1^2x_3^2x_4^2 + x_2^2x_3^2x_4^2)
\end{array}
\label{eqn:indicator-function-F3}
\end{equation}

There are several approaches to consider the indicator functions of
multi-level designs. In \cite{Pistone-Rogantin-2008b}, a complex coding
is proposed to generalize the arguments on two-level cases to multi-level
cases. For example, instead of $\{-1,0,1\}$ above, 
the three-level factor is coded as $\{1, w, w^2\}$, where $w =
\exp(2\pi\sqrt{-1}/3)$ in
\cite{Pistone-Rogantin-2008b}. The idea of the complex coding is based
on a theory of a harmonic analysis, where the indicator function is
viewed as a discrete Fourier transform.
Other approach is presented in
\cite{Cheng-Yw-2004} for the real coefficients field. However, it is
better if we can consider 
$\mathbb{Q}$, the field of rational numbers, as the coefficients field, 
because algebraic computations are conducted in $\mathbb{Q}$ (or finite
fields $\mathbb{Z}/p\mathbb{Z}$) for standard computational algebraic
software. 
Another approach is a concept of Hilbert basis presented in 
\cite{Carlini-Pistone-2009}, where the case of repeated treatments are
considered by considering counting functions instead of indicator
functions.
In this paper, 
we give generalization of the relations for two-level designs such as 
(\ref{eqn:relation-theta-2-level}) and
(\ref{eqn:example-constraints-2-level}) to general multi-level designs
for the rational coefficients field $\mathbb{Q}$,
and show how to relate the structure of the designs to the
structure of their indicator functions. 

The construction of this paper is as follows. In Section 
\ref{sec:definition-indicator-function}, we give
necessary definitions and theorems on the indicator functions. 
In Section 3, we give the structure of the indicator functions for
general fractional factorial designs. 
We also derive another
representation of the indicator functions, namely, contrast representation,
to reflect the orthogonality of the designs directly. 
In Section 4, we use these results to classify $2^3\times 3$ and
$2^4\times 3$ designs with given orthogonalities by a computational
algebraic software.

\section{The indicator functions of fractional factorial designs}
\label{sec:definition-indicator-function}
In this section, we give necessary materials on the indicator functions
of fractional factorial designs. The arguments are based on the theory
of interpolatory polynomial functions on designs, which is one of the
first applications of Gr\"obner basis theory to statistics 
introduced by \cite{Pistone-Wynn-1996}. 
See \cite{Pistone-Riccomagno-Wynn-2001} or Chapter 5 of 
\cite{Cox-Little-OShea-1992} for detail.

Let $x_1,\ldots,x_n$ be $n$ factors. Let $A_j \subset
\mathbb{Q}$ be a level set of a factor $x_j$ for $j = 1,\ldots,n$, 
where $\mathbb{Q}$ denotes the field of rational numbers.
We denote by $r_j = \#A_j$ the cardinality of $A_j$ and assume $r_j \geq
2$ for $j = 1,\ldots,n$.  
{\it A full factorial design} of the factors $x_1,\ldots,x_n$ is $D =
A_1\times \cdots\times A_n$. 
For later use, we introduce an {\it index set}
\[
 {\cal I} = \{(i_1,\ldots,i_n) \in [r_1]\times \cdots \times [r_n]\},
\]
where $[k] = \{1,2,\ldots,k\}$ for a positive integer $k$. We specify
each point of $D$ as $D = \{\Bd_{\Bi} \in \mathbb{Q}^n\ :\ \Bi \in {\cal
I}\}$. When we code $A_j = [r_j]$ for $j = 1,\ldots,n$, ${\cal I}$ coincides
with $D$ itself.  

A subset of $D$ is called {\it a fractional factorial design}. 
A fractional factorial design $F \subset D$ can be written as $F =
\{\Bd_{\Bi} \in D\ :\ \Bi \in {\cal I}'\}$ where ${\cal I}'$ is a subset of
${\cal I}$. 
Each design can be viewed as a finite subset of
$\mathbb{Q}^n$, i.e., as an algebraic variety, because each design can
be characterized as the set of the solutions of a system of polynomial
equations with rational coefficients. 
The {\it size} of a design is the cardinality of the design. 
We write the
size of a full factorial design $D$ as $m = \prod_{j = 1}^n r_j$ for
later use.

Let $\mathbb{Q}[x_1,\ldots,x_n]$ be the polynomial ring with
coefficients in $\mathbb{Q}$. For a design $F  \subset
\mathbb{Q}^n$, we denote by $I(F)$ the set of polynomials in
$\mathbb{Q}[x_1,\ldots,x_n]$ which are $0$ at every point of $F$, i.e., 
\[
 I(F) = \{f \in \mathbb{Q}[x_1,\ldots,x_n]\ :\ f(\Bd) = 0,\ \forall \Bd
 \in F\}.
\]
It is easy to prove that the set $I(F)$ is an ideal of
$\mathbb{Q}[x_1,\ldots,x_n]$. $I(F)$ is called {\it the design ideal} of
$F$. The design ideal introduced by \cite{Pistone-Wynn-1996} is a
fundamental tool to consider designs
algebraically. 
The design ideal is a radical ideal (Theorem 20 of
\cite{Pistone-Riccomagno-Wynn-2001}). The set of points $\Bd
\in \mathbb{Q}^n$ satisfying $f(\Bd) = 0$ for all $f \in I(F)$ is $F$ itself. 

For a full factorial design $D$, the design ideal $I(D)$ can be written
as
\[
 I(D) = \left<x_j^{r_j} - g_j,\ j = 1,\ldots,n \right>,
\]
where $g_j$ is a polynomial in $\mathbb{Q}[x_j]$ with the
degree less than $r_j$, $j = 1,\ldots,n$. 
Here $\left< \{f_i\} \right>$ means ``the ideal generated by
$\{f_i\}$''. 
In other words, 
the set
\begin{equation}
 G = \left\{x_j^{r_j} - g_j,\ j = 1,\ldots,n \right\}
\label{eqn:G-I(D)}
\end{equation}
is a generator of $I(D)$. In addition, $G$ is a 
reduced Gr\"obner basis of $I(D)$ for any monomial order. 
Note that the term $x_j^{r_j}$ is greater than the leading term of $g_j$
with respect to any monomial order.
We write the set of the monomials that are not divisible by the initial
monomials of $G$, $\{x_j^{r_j},\ j = 1,\ldots,n\}$, as
\[
 {\rm Est}(D) = \left\{\Bx^{\Ba} = \prod_{j = 1}^nx_j^{a_j}\ :\ \Ba \in
 L\right\},
\]
where 
\[
 L = \{\Ba = (a_1,\ldots,a_n) \in \mathbb{Z}_{\geq 0}^n\ :\ 0 \leq a_j
 \leq r_j-1,\ j = 1,\ldots,n\}
\]
and $\mathbb{Z}_{\geq 0}$ is the set of nonnegative integers. 
Note that the cardinality of $L$ is $m$. 
From $D = \{\Bd_{\Bi} \in \mathbb{Q}^n\ :\ \Bi \in {\cal I}\}$ and $L$, we
define {\it a model matrix} by 
\[
 X = \left[\Bd_{\Bi}^{\Ba}\right]_{\Bi \in {\cal I}; \Ba \in L}, 
\]
where $\Bd_{\Bi}^{\Ba} = \prod_{j = 1}^nd_{\Bi j}^{a_j}$ and $d_{\Bi j}$
is the level of the factor $j$ in the experimental run indexed by $\Bi
\in {\cal I}$. Note that $X$ is called a design matrix in Definition 26
of \cite{Pistone-Riccomagno-Wynn-2001}.
By ordering the elements of ${\cal I}$ and $L$, $X$ is an $m \times m$
matrix, and is nonsingular (Theorem 26 of
\cite{Pistone-Riccomagno-Wynn-2001}).

The quotient of $\mathbb{Q}[x_1,\ldots,x_n]$ modulo the  
design ideal $I(D) \in \mathbb{Q}[x_1,\ldots,x_n]$ is defined by 
\[
 \mathbb{Q}[x_1,\ldots,x_n]/I(D) = \{[f]\ :\ f \in
 \mathbb{Q}[x_1,\ldots,x_n]\},
\]
where we define $[f] = \{g \in \mathbb{Q}[x_1,\ldots,x_n]\ {\rm such\ 
that}\ f - g \in I(D)\}$. In the terminology of the designs of
experiments, two polynomial models $f$ and $g$ are {\it confounded} on $D$ if
and only if $f - g \in I(D)$. 
Therefore each element $[f] \in  \mathbb{Q}[x_1,\ldots,x_n]/I(D)$ is
the set of the polynomials $g \in \mathbb{Q}[x_1,\ldots,x_n]$ that is
confounded to $f$ on $D$. 
An important fact is that 
${\rm Est}(D)$ is a basis of $\mathbb{Q}[x_1,\ldots,x_n]/I(D)$ as a
$\mathbb{Q}$-vector space. See Theorem 15 of
\cite{Pistone-Riccomagno-Wynn-2001} for detail.

Suppose we have a $\mathbb{Q}$-valued response (or, observations) $\By =
(y(\Bi))_{\Bi \in {\cal I}}$ on $D = \{\Bd_i\in \mathbb{Q}^n:\ \Bi \in
{\cal I}\}$. 
Note that each
polynomial $f \in \mathbb{Q}[x_1,\ldots,x_n]$ can be viewed as 
a response function on $D$, i.e., 
$f \in \mathbb{Q}^D$ where we
denote $\mathbb{Q}^D$
by the vector space of functions from $D$ to $\mathbb{Q}$. 
The {\it interpolatory polynomial function} for $\By$ is a polynomial $f
\in \mathbb{Q}[x_1,\ldots,x_n]$ satisfying $f(\Bd_{\Bi}) = y(\Bi),\ \Bi
\in {\cal I}$. From the fact that ${\rm Est}(D)$ is a basis of
$\mathbb{Q}[x_1,\ldots,x_n]/I(D)$, the interpolatory polynomial function
for $\By$ is written uniquely as 
\begin{equation}
 f(x_1,\ldots,x_n) = \sum_{\Ba \in L}\theta_{\Ba}\Bx^{\Ba}, 
\label{eqn:poly-inter-polatory}
\end{equation}
where an $m \times 1$ column vector $\Btheta = (\theta_{\Ba})_{\Ba \in
L}$ is given by $\Btheta = X^{-1}\By$ for an $m\times 1$ column vector
$\By$. See Theorem 26 of \cite{Pistone-Riccomagno-Wynn-2001} for detail.

Now we introduce an indicator function. 
\begin{definition}[\cite{Fontana-Pistone-Rogantin-2000}]
Let $F \subset D$ be a fractional factorial design. 
The indicator function of $F$ is a response function $f$ on $D$
 satisfying
\[
 f(\Bd) = \left\{\begin{array}{ll}
1, & \mbox{if}\ \Bd \in F,\\
0, & \mbox{if}\ \Bd \in D \setminus F.
\end{array}
\right.
\]
\end{definition}
By the definition, the indicator function is constructed as follows. Write a 
fractional factorial design $F \subset D$ as $F = \{\Bd_{\Bi} \in D\ :\
\Bi \in {\cal I}'\}$ for a subset ${\cal I}' \subset {\cal I}$. Then the
indicator function of $F$ is the interpolatory polynomial function for a
response $\By = (y(\Bi))_{\Bi \in {\cal I}}$, where
\begin{equation}
 y(\Bi) = \left\{\begin{array}{ll}
1, & \mbox{if}\ \Bi \in {\cal I}'\\
0, & \mbox{if}\ \Bi \in {\cal I} \setminus {\cal I}'.
\end{array}
\right.
\label{eqn:y-subset-design}
\end{equation}
From the uniqueness of the interpolarory polynomial function mentioned
above, the representation of the indicator function is unique.

\begin{example}\label{example:2x2x3}
Consider designs of $3$ factors $x_1,x_2,x_3$, where $x_1, x_2$ are
 two-level factors and $x_3$ is a three-level factor. 
We code the levels of each factor as
\[
 A_1 = A_2 = \{-1,1\},\ A_3 = \{-1,0,1\}.
\]
Therefore $r_1 = r_2 = 2, r_3 = 3$, and 
the full factorial design $D = A_1\times A_2\times A_3$ has $m = 12$ points.
The index set is ${\cal I} = \{1,2\}\times \{1,2\}\times \{1,2,3\}$. 
The full factorial design $D$ is written as $D = \{\Bd_{\Bi}\ :\ \Bi \in
 {\cal I}\}$, where
\[
 \begin{array}{lll}
\Bd_{(1,1,1)} = (-1,-1,-1), & \Bd_{(1,1,2)} = (-1,-1,0), &
 \Bd_{(1,1,3)} = (-1,-1,1), \\
\Bd_{(1,2,1)} = (-1,1,-1), & \Bd_{(1,2,2)} = (-1,1,0), &
 \Bd_{(1,2,3)} = (-1,1,1), \\
\Bd_{(2,1,1)} = (1,-1,-1), & \Bd_{(2,1,2)} = (1,-1,0), &
 \Bd_{(2,1,3)} = (1,-1,1), \\
\Bd_{(2,2,1)} = (1,1,-1), & \Bd_{(2,2,2)} = (1,1,0), &
 \Bd_{(2,2,3)} = (1,1,1).
 \end{array}
\]
The design ideal of $D$ is written as
\[
 I(D) = \left<x_1^2-1, x_2^2-1, x_3^3 - x_3\right> \subset
 \mathbb{Q}[x_1,x_2,x_3], 
\]
and $G = \{x_1^2-1, x_2^2-1, x_3^3 - x_3\}$ is a reduced Gr\"obner basis
 of $I(D)$ for any monomial order.
Therefore we have
\[
 {\rm Est}(D) =
 \{1,x_1,x_2,x_3,x_3^2,x_1x_2,x_1x_3,x_2x_3,x_1x_2x_3,x_1x_3^2,x_2x_3^2,
x_1x_2x_3^2 \}.
\]
Note that there are $m$ monomials in ${\rm Est}(D)$. Corresponding $L$
 is given by
\[\begin{array}{rcl}
 L & = & \{\ (0,0,0), (1,0,0), (0,1,0), (0,0,1), (0,0,2), (1,1,0),\\
& & \hspace*{4mm}(1,0,1), (0,1,1), (1,1,1), (1,0,2), (0,1,2), (1,1,2)\ \}.
\end{array}\]
The model matrix $X$ is given in Figure
 \ref{fig:model-matrix-2x2x3}. 
Note that, differently than in the usual ANOVA decomposition, the
 columns corresponding to the interactions are not orthogonal to the
 columns corresponding to the simple factors.
Here, and hereafter, we write each element
 of $L$ and ${\cal I}$ by omitting commas as $(a_1\cdots a_n)$ or
 $(i_1\cdots i_n)$ instead
 of $(a_1,\cdots,a_n)$ or $(i_1,\cdots,i_n)$ for simplicity.
\begin{figure*}[htbp]
\[{\small
\begin{array}{c|rrrrrrrrrrrr|}
\multicolumn{1}{c}{{\cal I}\backslash L} & 000 & 100 & 010 & 001 & 002 & 110 &
101 & 011 & 111 & 102 & 012 & \multicolumn{1}{r}{112}\\ \cline{2-13}
111 &1&-1&-1&-1&1& 1& 1& 1&-1&-1&-1& 1\\ 
112 &1&-1&-1& 0&0& 1& 0& 0& 0& 0& 0& 0\\ 
113 &1&-1&-1& 1&1& 1&-1&-1& 1&-1&-1& 1\\ 
121 &1&-1& 1&-1&1&-1& 1&-1& 1&-1& 1&-1\\ 
122 &1&-1& 1& 0&0&-1& 0& 0& 0& 0& 0& 0\\ 
123 &1&-1& 1& 1&1&-1&-1& 1&-1&-1& 1&-1\\
211 &1& 1&-1&-1&1&-1&-1& 1& 1& 1&-1&-1\\ 
212 &1& 1&-1& 0&0&-1& 0& 0& 0& 0& 0& 0\\ 
213 &1& 1&-1& 1&1&-1& 1&-1&-1& 1&-1&-1\\ 
221 &1& 1& 1&-1&1& 1&-1&-1&-1& 1& 1& 1\\ 
222 &1& 1& 1& 0&0& 1& 0& 0& 0& 0& 0& 0\\ 
223 &1& 1& 1& 1&1& 1& 1& 1& 1& 1& 1& 1\\ \cline{2-13}
\end{array}
}\]
\caption{The model matrix of a full factorial design $D =
 \{-1,1\}^2\times \{-1,0,1\}$.}
\label{fig:model-matrix-2x2x3}
\end{figure*}

Now consider a fractional factorial design $F = \{\Bd_{\Bi}\ :\ \Bi \in
 {\cal I}'\}\subset D$, where ${\cal I}' =
 \{(111),(122),(213),(223)\}$. 
The indicator function of $F$ is constructed as the  
interpolatory
 polynomial function (\ref{eqn:poly-inter-polatory}) for a response $\By
 = (y(\Bi))_{\Bi \in {\cal I}}$
 satisfying (\ref{eqn:y-subset-design}). For this ${\cal I}'$, 
$\By = (1,0,0,0,1,0,0,0,1,0,0,1)^{T}$ and 
we have
\[
 \Btheta  =  X^{-1}\By 
 =  \frac{1}{8}(2,-2,2,1,1,-2,3,1,-1,3,-3,3)^{T},
\]
where ${}^{T}$ is a transpose. Therefore the indicator function of $F$
 is
\[\begin{array}{l}
 f(x_1,x_2,x_3) = \displaystyle \frac{1}{4} - \frac{1}{4}(x_1 - x_2) +
 \frac{1}{8}(x_3 + x_3^2) - \frac{1}{4}x_1x_2 + \frac{3}{8}x_1x_3 +
 \frac{1}{8}x_2x_3\vspace*{2mm}\\
\multicolumn{1}{r}{\displaystyle {} - \frac{1}{8}x_1x_2x_3 +
 \frac{3}{8}(x_1x_3^2 -
 x_2x_3^2 + x_1x_2x_3^2).}
\end{array}\]
\hspace*{\fill}$\Box$
\end{example}

\section{Characterization of orthogonal fractional factorial designs by
 indicator functions}
Now we consider relations between the design and its indicator
function. We start with the generalization of the relation
(\ref{eqn:relation-theta-2-level}) for two-level designs to
multi-level designs. 

A polynomial $f \in \mathbb{Q}[x_1,\ldots,x_n]$ is an indicator function
of some fractional factorial design $F \subset D = \{\Bd_{\Bi} \in
\mathbb{Q}^n\ :\ \Bi \in {\cal I}\}$ if and only if $f^2 - f \in I(D)$,
i.e., $f$ and $f^2$ are in the same equivalence class of
$\mathbb{Q}[x_1,\ldots,x_n]/I(D)$. Therefore, suppose $f$ represented as 
(\ref{eqn:poly-inter-polatory}) 
is an indicator function of some
fractional factorial design, we have
\[\begin{array}{rcl}
\displaystyle \sum_{\Ba \in L}\theta_{\Ba}\Bx^{\Ba} 
& = & \displaystyle
 \left(\sum_{\Ba \in
 L}\theta_{\Ba}\Bx^{\Ba}\right)^2\ \ \ {\rm mod}\ I(D)\\
& = & \displaystyle
\sum_{\Ba_1 \in
 L}\sum_{\Ba_2 \in L}\theta_{\Ba_1}\theta_{\Ba_2}\Bx^{\Ba_1 + \Ba_2}\ \
 \ {\rm mod}\ I(D).
\end{array}
\]
Here, write the standard form of 
$\displaystyle
\sum_{\Ba_1 \in
 L}\sum_{\Ba_2 \in L}\theta_{\Ba_1}\theta_{\Ba_2}\Bx^{\Ba_1 + \Ba_2}$
 with respect to $G$ as 
\begin{equation}
 r = \sum_{\Ba\in L}\mu_{\Ba}\Bx^{\Ba}.
\label{eqn:rhs-r-mu}
\end{equation}
In other words, $r$ is a unique remainder when we divide 
$\displaystyle
\sum_{\Ba_1 \in
 L}\sum_{\Ba_2 \in L}\theta_{\Ba_1}\theta_{\Ba_2}\Bx^{\Ba_1 + \Ba_2}$ by
 $G$, the reduced Gr\"ober basis of $I(D)$. Then we have the following result.
\begin{proposition}[Generalization of Proposition 3.7 of
 \cite{Fontana-Pistone-Rogantin-2000}]
A polynomial $f$ represented as 
(\ref{eqn:poly-inter-polatory}) is an indicator function of some
 fractional factorial design if and only if a system of algebraic
 equations
\begin{equation}
 \theta_{\Ba} = \mu_{\Ba},\ \ \ \Ba \in L
\label{eq:prop-theta-mu-eq}
\end{equation}
holds, where $\mu_{\Ba}$ is given by (\ref{eqn:rhs-r-mu}).
\end{proposition}

\paragraph*{Proof.} From the division algorithm and the property of the 
Gr\"obner basis. 
See Chapter 2 of \cite{Cox-Little-OShea-1992}.\hspace*{\fill}$\Box$

\bigskip

The meaning of the relation (\ref{eq:prop-theta-mu-eq}) is explained as
follows. From the theory of the interpolatory polynomial function on
$D$, each $\mathbb{Q}$-valued function on $D$ is
represented (on $D$) by a polynomial with the monomials in ${\rm Est}(D)$.
This applies to the indicator function $f$ of a fraction $F \subset
D$. In this case, the reduced Gr\"obner basis of $I(D)$ is simply the
list of the univariate monic polynomials (\ref{eqn:G-I(D)}) defining
the levels of each factor, and the remainder $r$ is derived by
substitution as it is. 

\begin{example}[Continuation of Example \ref{example:2x2x3}]
\label{example:2x2x3-2}
Consider $2\times 2\times 3$ designs. When we code the levels as $A_1 =
 A_2 = \{-1,1\}, A_3 = \{-1,0,1\}$, the relation
 (\ref{eq:prop-theta-mu-eq}) is as follows.
\[
 \begin{array}{rcl}
\theta_{000} & = & \theta_{100}^2+\theta_{010}^2+\theta_{000}^2+\theta_{110}^2\\
\theta_{100} & = & 2\theta_{100}\theta_{000}+2\theta_{010}\theta_{110}\\
\theta_{010} & = & 2\theta_{100}\theta_{110}+2\theta_{010}\theta_{000}\\
\theta_{001} & = & 2\theta_{100}\theta_{101} + 2\theta_{010}\theta_{011}
 + 2\theta_{001}\theta_{002} + 2\theta_{001}\theta_{000} +
 2\theta_{110}\theta_{111} +
 2\theta_{101}\theta_{102}+2\theta_{011}\theta_{012} +
 2\theta_{111}\theta_{112}\\
\theta_{002} & = & 2\theta_{100}\theta_{102} + 2\theta_{010}\theta_{012}
 + \theta_{001}^2+\theta_{002}^2 + 2\theta_{002}\theta_{000} +
 2\theta_{110}\theta_{112} + \theta_{101}^2+\theta_{011}^2 +
 \theta_{111}^2 + \theta_{102}^2 + \theta_{012}^2\\
& &  + \theta_{112}^2\\
\theta_{110} & = & 2\theta_{100}\theta_{010} + 2\theta_{000}\theta_{110}\\
\theta_{101} & = & 2\theta_{100}\theta_{001} + 2\theta_{010}\theta_{111}
 + 2\theta_{001}\theta_{102} + 2\theta_{002}\theta_{101} +
 2\theta_{000}\theta_{101} + 2\theta_{110}\theta_{011} +
 2\theta_{011}\theta_{112} + 2\theta_{111}\theta_{012}\\
\theta_{011} & = & 2\theta_{100}\theta_{111} + 2\theta_{010}\theta_{001}
 + 2\theta_{001}\theta_{012} + 2\theta_{002}\theta_{011} +
 2\theta_{000}\theta_{011} + 2\theta_{110}\theta_{101} +
 2\theta_{101}\theta_{112} + 2\theta_{111}\theta_{102}\\
\theta_{111} & = & 2\theta_{100}\theta_{011} + 2\theta_{010}\theta_{101}
 + 2\theta_{001}\theta_{110} + 2\theta_{001}\theta_{112} +
 2\theta_{002}\theta_{111} + 2\theta_{000}\theta_{111} +
 2\theta_{101}\theta_{012} + 2\theta_{011}\theta_{102}\\
\theta_{102} & = & 2\theta_{100}\theta_{002} + 2\theta_{010}\theta_{112}
 + 2\theta_{001}\theta_{101} + 2\theta_{002}\theta_{102} +
 2\theta_{000}\theta_{102} + 2\theta_{110}\theta_{012} +
 2\theta_{011}\theta_{111} + 2\theta_{012}\theta_{112}\\
\theta_{012} & = & 2\theta_{100}\theta_{112} + 2\theta_{010}\theta_{002}
 + 2\theta_{001}\theta_{011} + 2\theta_{002}\theta_{012} +
 2\theta_{000}\theta_{012} + 2\theta_{110}\theta_{102} +
 2\theta_{101}\theta_{111} + 2\theta_{102}\theta_{112}\\
\theta_{112} & = & 2\theta_{100}\theta_{012} + 2\theta_{010}\theta_{102}
 + 2\theta_{001}\theta_{111} + 2\theta_{002}\theta_{110} +
 2\theta_{002}\theta_{112} + 2\theta_{000}\theta_{112} +
 2\theta_{101}\theta_{011} + 2\theta_{102}\theta_{012}
 \end{array}
\]
These results are easily obtained by standard algebraic softwares, such as 
Macaulay2 (\cite{Macaulay2}). 
Each solution of the above system of polynomial equations corresponds to
 the coefficients of the indicator function for each fractional
 factorial designs. 

If we change the level codings, the relation (\ref{eq:prop-theta-mu-eq})
 changes. For example, when we code the levels as $A_1 = A_2 = \{0,1\},
 A_3 = \{0,1,2\}$, the relation (\ref{eq:prop-theta-mu-eq}) is as follows.
\[
\begin{array}{rcl}
\theta_{000} & = & \theta_{000}^2\\
\theta_{100} & = & \theta_{100}^2+2\theta_{100}\theta_{000}\\
\theta_{010} & = & \theta_{010}^2+2\theta_{010}\theta_{000}\\
\theta_{001} & = &
 -4\theta_{001}\theta_{002}-6\theta_{002}^2+2\theta_{001}\theta_{000}\\
\theta_{002} & = & \theta_{001}^2 + 6\theta_{001}\theta_{002} +
 7\theta_{002}^2 +
 2\theta_{002}\theta_{000}\\
\theta_{110} & = & 2\theta_{100}\theta_{010} + 2\theta_{100}\theta_{110} +
 2\theta_{010}\theta_{110}
 + 2\theta_{000}\theta_{110} + \theta_{110}^2\\
\theta_{101} & = & 2\theta_{100}\theta_{001} + 2\theta_{100}\theta_{101} -
 4\theta_{002}\theta_{101} + 2\theta_{000}\theta_{101} -
 4\theta_{001}\theta_{102} - 12\theta_{002}\theta_{102} -
 4\theta_{101}\theta_{102} - 6\theta_{102}^2\\
\theta_{011} & = & 2\theta_{010}\theta_{001} + 2\theta_{010}\theta_{011} -
 4\theta_{002}\theta_{011} + 2\theta_{000}\theta_{011} -
 4\theta_{001}\theta_{012} - 12\theta_{002}\theta_{012} -
 4\theta_{011}\theta_{012} - 6\theta_{012}^2\\
\theta_{111} & = & 2\theta_{001}\theta_{110} + 2\theta_{010}\theta_{101} +
 2\theta_{110}\theta_{101} + 2\theta_{100}\theta_{011} +
 2\theta_{110}\theta_{011} + 2\theta_{100}\theta_{111} +
 2\theta_{010}\theta_{111} - 4\theta_{002}\theta_{111}\\
& &  {} + 2\theta_{000}\theta_{111} + 2\theta_{110}\theta_{111} -
 4\theta_{011}\theta_{102} - 4\theta_{111}\theta_{102} -
 4\theta_{101}\theta_{012} - 4\theta_{111}\theta_{012} -
 12\theta_{102}\theta_{012}\\
& & {} - 4\theta_{001}\theta_{112} -
 12\theta_{002}\theta_{112} - 4\theta_{101}\theta_{112} -
 4\theta_{011}\theta_{112} - 4\theta_{111}\theta_{112} -
 12\theta_{102}\theta_{112} - 12\theta_{012}\theta_{112} -
 6\theta_{112}^2\\
\theta_{102} & = & 2\theta_{100}\theta_{002} + 2\theta_{001}\theta_{101} +
 6\theta_{002}\theta_{101} + \theta_{101}^2 + 2\theta_{100}\theta_{102}
 + 6\theta_{001}\theta_{102} + 14\theta_{002}\theta_{102} +
 2\theta_{000}\theta_{102}\\
& &  + 6\theta_{101}\theta_{102} + 7\theta_{102}^2\\
\theta_{012} & = & 2\theta_{010}\theta_{002} + 2\theta_{001}\theta_{011} +
 6\theta_{002}\theta_{011} + \theta_{011}^2 + 2\theta_{010}\theta_{012}
 + 6\theta_{001}\theta_{012} + 14\theta_{002}\theta_{012} +
 2\theta_{000}\theta_{012}\\
& & {}+ 6\theta_{011}\theta_{012} +
 7\theta_{012}^2\\
\theta_{112} & = & 2\theta_{002}\theta_{110} + 2\theta_{101}\theta_{011} +
 2\theta_{001}\theta_{111} + 6\theta_{002}\theta_{111} +
 2\theta_{101}\theta_{111} + 2\theta_{011}\theta_{111} + \theta_{111}^2
 + 2\theta_{010}\theta_{102}\\
& & + 2\theta_{110}\theta_{102} +
 6\theta_{011}\theta_{102} + 6\theta_{111}\theta_{102} +
 2\theta_{100}\theta_{012} + 2\theta_{110}\theta_{012} +
 6\theta_{101}\theta_{012} + 6\theta_{111}\theta_{012}\\
& & {}+ 14\theta_{102}\theta_{012} + 2\theta_{100}\theta_{112} +
 2\theta_{010}\theta_{112} + 6\theta_{001}\theta_{112} +
 14\theta_{002}\theta_{112} + 2\theta_{000}\theta_{112} +
 2\theta_{110}\theta_{112}\\
& & {}+ 6\theta_{101}\theta_{112} +
 6\theta_{011}\theta_{112} + 6\theta_{111}\theta_{112} +
 14\theta_{102}\theta_{112} + 14\theta_{012}\theta_{112} +
 7\theta_{112}^2\\
\end{array}
\]
In actual applications, there are cases 
where the level coding has not essential meaning, such as for the
 designs of qualitative factors. However, for our purpose of solving a
 system of polynomial equations using computational algebraic software,
an appropriate level coding is important in view of computational
 time. In the author's experiences, it is better to code $\{-1,1\}$ rather than
 $\{0,1\}$ for two-level factor, and $\{-1,0,1\}$ rather than
 $\{0,1,2\}$ for three-level factor. 
We consider this point in Section
 \ref{sec:computation}.
\hspace*{\fill}$\Box$
\end{example}

\bigskip

As we see in Example \ref{example:2x2x3-2}, the relation 
(\ref{eq:prop-theta-mu-eq}) is very complicated compared to the relation
for two-level cases (\ref{eqn:relation-theta-2-level}). 
Among the various characterizations of the
coefficients of the indicator
functions of two-level designs given in
\cite{Fontana-Pistone-Rogantin-2000}, the relation of the indicator
functions of complementary designs can be generalized to multi-level cases
as follows.
\begin{proposition}[Generalization of Corollary 3.5 of
 \cite{Fontana-Pistone-Rogantin-2000}]
If $F$ and $\bar{F}$ are complementary fractions and $\Btheta =
 (\theta_{\Ba})_{\Ba \in L}$ and $\bar{\Btheta} =
 (\bar{\theta}_{\Ba})_{\Ba \in L}$ are the coefficients of the
 corresponding 
indicator functions given by
 (\ref{eqn:poly-inter-polatory}) respectively, then 
\[
 \theta_{0\cdots 0} = 1 - \bar{\theta}_{0\cdots 0}\ \ \mbox{and}\ \
 \theta_{\Ba} = -\bar{\theta}_{\Ba},\ \forall\ \Ba \neq (0,\ldots,0).
\]
\end{proposition}
\paragraph*{Proof.}\ Write the model matrix $X$ as
\[
 X = \left[\Bone_m\ |\ \Bs_2\ |\ \cdots\ |\ \Bs_m\right]
\]
and write $X^{-1}\Bone_m = (c_1,\ldots,c_m)^{T}$. Then we have $c_1 =
1,\ c_2 = \cdots = c_m = 0$ from the non-singularity of $X$ in the
relation
\[
 \Bone_m = c_1\Bone_m + c_2\Bs_2 + \cdots + c_m\Bs_m.
\]
Therefore for the responses $\By, \bar{\By} \in \{0,1\}^m$ such that
$\By + \bar{\By} = \Bone_m$, we have
\[
 \Btheta + \bar{\Btheta} = X^{-1}(\By + \bar{\By}) =
 (1,0,\ldots,0)^T.
\]
\hspace*{\fill}$\Box$

\bigskip

Next consider structure of the indicator functions of designs with given
characteristic. 
The idea is to express the structure of the indicator functions as
additional constraints to the system of polynomial equations
(\ref{eq:prop-theta-mu-eq}) to classify designs with given
characteristic. The additional constraints are derived as
follows. Recall that the coefficients vector $\Btheta$ is given by
$\Btheta = X^{-1}\By$ in (\ref{eqn:poly-inter-polatory}). Here, treat 
$\By = (y(\Bi))_{\Bi \in {\cal I}}$ as a vector of $\{0,1\}^m$ in
(\ref{eqn:y-subset-design}) and express the characteristic of designs as
\begin{equation}
 \Bc^{T}\By = s,\ s \in \mathbb{Q}
\label{eqn:add-const-cy}
\end{equation}
for a constant column vector $\Bc \in \mathbb{Q}^m$.
For example, $\By \in \{0,1\}^m$ corresponding to designs with the size $s$
satisfies the constraint 
\[
 \Bone_m^{T}\By = s,
\] 
where $\Bone_m = (1,\ldots,1)^T$ is an $m\times 1$ column vector of the
elements $1$'s. Equireplicated designs or orthogonal designs can be
expressed by 
\[
 \Bc^T\By = 0
\] 
for some contrast vectors $\Bc$. 

Based on the above idea, 
we define {\it a contrast matrix}. 
For each subset $J \subset [n]$, we define 
${\cal I}_J = \prod_{j \in J}[r_j] \subset {\cal I}$ and its cardinality
by $m_J = \prod_{j \in J}r_j$.
We also define $\Bi_J$ by the restriction of $\Bi = (i_1,\ldots,i_n) \in
{\cal I}$ to the index of ${\cal
I}_J$. For example of $n = 4$ and $J = \{1,2,4\}$, we have 
$\Bi_J = (i_1,i_2,i_4)$. 

\begin{definition}\label{def:contrast-matrix}
The contrast matrix $C$ is an $m\times m$
 matrix of the form
\[
 C^{T} =\left[\Bone_m\ |\ C_1^T\ |\ C_2^T\ |\ \cdots\ |\ C_n^T\right], 
\]
where $C_k$ is a $v_k \times m$ matrix where
\[
 v_k = \sum_{J \subset [n], \# J = k}\left(
\prod_{j \in J}(r_j-1)
\right).
\]
The set of $m\times 1$ column vectors of $C_k^T$ is 
\[
\left\{\Bc_{J(\tilde{\Bi})} = \{c_{J(\tilde{\Bi})}(\Bi)\}_{\Bi \in {\cal I}}\ :\ J
 \subset [n], \#J = k, \tilde{\Bi} \in \prod_{j \in J}[r_j-1]\right\},
\]
where
\[
 c_{J(\tilde{i})}(\Bi) = \left\{\begin{array}{rl}
1, &  \Bi_J = 1,\\
-1, & \Bi_J = \tilde{i} + 1,\\
0, & \mbox{otherwise}
\end{array}
\right.
\]
for $\#J = 1$, and
\[
 c_{J(\tilde{\Bi})}(\Bi) = \left\{\begin{array}{rl}
1, & \Bi_J = (\tilde{i}_1,\ldots,\tilde{i}_{k-1},1),\\
-1, & \Bi_J = (\tilde{i}_1,\ldots,\tilde{i}_{k-1},\tilde{i}_k + 1),\\
0, & \mbox{otherwise}
\end{array}
\right.
\]
for $\#J \geq 2$. 
\end{definition}

\bigskip

Note that the contrast matrix $C$ is constructed only from ${\cal I}$. 
In other words, $C$ is uniquely determined from
$r_1,\ldots,r_n$.  Especially, $C$ does not depend on the level coding.

\bigskip

\begin{example}[Continuation of Example \ref{example:2x2x3-2}]
\label{example:2x2x3-3}
For $2 \times 2\times 3$ designs, the contrast matrix $C$ is given in
 Figure \ref{fig:contrast-matrix-2x2x3}.
\begin{figure*}[htbp]
\[{\small
\begin{array}{c|rrrrrrrrrrrr|}
\multicolumn{1}{c}{J(\tilde{\Bi})\backslash {\cal I}} & 111 & 112 & 113
 & 121 & 122 & 123 &
211 & 212 & 213 & 221 & 222 & \multicolumn{1}{r}{223}\\ \cline{2-13}
{\rm Const.} &1& 1& 1& 1& 1& 1& 1& 1& 1& 1& 1& 1\\ 
1(1)     &1& 1& 1& 1& 1& 1&-1&-1&-1&-1&-1&-1\\ 
2(1)     &1& 1& 1&-1&-1&-1& 1& 1& 1&-1&-1&-1\\ 
3(1)     &1&-1& 0& 1&-1& 0& 1&-1& 0& 1&-1& 0\\ 
3(2)     &1& 0&-1& 1& 0&-1& 1& 0&-1& 1& 0&-1\\ 
12(11)   &1& 1& 1&-1&-1&-1& 0& 0& 0& 0& 0& 0\\
13(11)   &1&-1& 0& 1&-1& 0& 0& 0& 0& 0& 0& 0\\ 
13(12)   &1& 0&-1& 1& 0&-1& 0& 0& 0& 0& 0& 0\\ 
23(11)   &1&-1& 0& 0& 0& 0& 1&-1& 0& 0& 0& 0\\ 
23(12)   &1& 0&-1& 0& 0& 0& 1& 0&-1& 0& 0& 0\\ 
123(111) &1&-1& 0& 0& 0& 0& 0& 0& 0& 0& 0& 0\\ 
123(112) &1& 0&-1& 0& 0& 0& 0& 0& 0& 0& 0& 0\\ \cline{2-13}
\end{array}
}\]
\caption{The contrast matrix of $2\times 2\times 3$ designs.}
\label{fig:contrast-matrix-2x2x3}
\end{figure*}
\hspace*{\fill}$\Box$
\end{example}

\bigskip

The contrast matrix $C$ given in Definition \ref{def:contrast-matrix} 
relates to the theory of contingency tables. Suppose the response $\By$ is a
vector of nonnegative integers, then we can treat $\By$ as a frequency
of contingency table 
$\By = \{y_{\Bi}\ :\ \Bi \in {\cal I}\}$ with the set of cells ${\cal
I}$. In this case, we see that the condition
\[
 \Bone_m^T \By = s,\ C_{\ell}\By = \Bzero_{v_{\ell}},\ \ell = 1,\ldots,t,
\]
means equal $\ell$-dimensional marginal totals for $\ell = 1,\ldots,t$. 
The definition of the contrast matrix and 
the following theorem (Theorem \ref{prop:orthogonal}) are based on this
connection. 
As another relation, $C$ is a
 configuration matrix in the theory of toric ideals. See Section1.5.3 of
 \cite{dojo-en}. 


By the contrast matrix, we specify the size and the orthogonality of the
designs as follows. We call a design $F \subset D$ is {\it orthogonal of
strength} $t$\ ($t \leq n$), if for any $t$ factors, all possible
combinations of levels appear equally often in $F$. This definition is
from the theory of orthogonal arrays. See Chapter 7 of
\cite{Wu-Hamada-2009}, for example. In particular, an orthogonal design
of strength $n$ is a full factorial design.
\begin{theorem}\label{prop:orthogonal}
If $\By \in \{0,1\}^m$ is a response on $D$ given by
 (\ref{eqn:y-subset-design}), the fractional factorial
 design $F = \{\Bd_{\Bi} \in D\ :\ \Bi \in {\cal I}'\}$ is 
size $s$ and orthogonal of strength $t$ if and only if
\begin{equation}
\left\{
\begin{array}{l}
 \Bone_m^T\By = s,\\
 C_k\By = \Bzero_{v_k},\ \ k = 1,\ldots,t,
\end{array}
\right.
\label{eqn:cond-for-orthogonal-y}
\end{equation}
where $\Bzero_{\ell} = (0,\ldots,0)^T$ is an $\ell\times 1$ column
 vector of the elements $0$'s, and $s$ is a common multiple of
 $\displaystyle\left\{\prod_{j \in
 J}r_j\ :\ \# J = t\right\}$.
\end{theorem}

To prove Theorem \ref{prop:orthogonal}, 
we define $J$-marginal vector of $\By =
(y(\Bi))_{\Bi \in {\cal I}}$ by $\By_J = (y_J(\Bi_J))_{\Bi_J \in {\cal
I}_J}$, where
\[
 y_J(\Bi_J) = \sum_{\Bi_{J^c} \in {\cal I}_{J^c}}y(\Bi_J,\Bi_{J^c}).
\]
Note that $J^c$ denotes the complement of $J$, and in $y(\Bi_J,
\Bi_{J^c})$, for notational simplicity, the indices in ${\cal I}_J$ are
collected to the left. Also we are writing $y(\Bi_J, \Bi_{J^c})$ instead
of $y((\Bi_J,\Bi_{J^c}))$.
For example of $n = 3$, we have $\By_{\{2\}} =
(y_{\{2\}(\Bi_{\{2\}})})_{\Bi_{\{2\}} \in {\cal I}_{\{2\}}} =
(y_{\{2\}}(i_2))_{i_2 \in [r_2]}$, where
\[
 y_{\{2\}}(i_2) = \sum_{i_1 \in [r_1]}\sum_{i_3 \in
 [r_3]}y(i_1,i_2,i_3), 
\]
and $\By_{\{1,3\}} = (y_{\{1,3\}}(\Bi_{\{1,3\}}))_{\Bi_{\{1,3\}} \in
{\cal I}_{\{1,3\}}} = (y_{\{1,3\}}(i_1,i_3))_{(i_1,i_3) \in
[r_1]\times[r_3]}$, where
\[
 y_{\{1,3\}}(i_1,i_3) = \sum_{i_2 \in [r_2]}y(i_1,i_2,i_3),
\]
and so on. The concept of the $J$-marginal vector is from the theory of
contingency tables. For detail, see Chapter 4.2 of \cite{Lauritzen}, ``Basic
facts and concepts'' of Contingency tables, for example. 

\paragraph*{Proof of Theorem \ref{prop:orthogonal}.}\ 
Suppose that the size $s$ is a common multiple of 
 $\displaystyle\left\{\prod_{j \in
 J}r_j\ :\ \# J = t\right\}$. 
Using the $J$-marginal vector $\By_J$, 
$F$ is size $s$ and orthogonal of
strength $t$ if and only if
\begin{equation}
\left\{
\begin{array}{l}
\Bone_m^T \By = s,\\
\displaystyle y_{J}(\Bi_J) = \frac{s}{m_J},\ \forall\ \Bi_J \in
{\cal I}_J\ \ \mbox{for all}\ J \subset [n]\ \mbox{with}\ \# J \leq t.
\end{array}
\right.
\label{eqn:cond-for-y_J-orthogonal}
\end{equation}
The relation $(\ref{eqn:cond-for-y_J-orthogonal}) \Rightarrow
(\ref{eqn:cond-for-orthogonal-y})$ is straightforward, because the
relation $C_k \By = \Bzero_{v_k}$ is equivalent to
\begin{equation}
 y_J(i_1,\ldots,i_{k-1},1) = 
 y_J(i_1,\ldots,i_{k-1},2) = \cdots = 
 y_J(i_1,\ldots,i_{k-1},r_k),\ \ \forall\  
(i_1,\ldots,i_{k-1}) \in \displaystyle \prod_{j = 1}^{k-1}[r_j - 1]
\label{eqn:y_J-edge}
\end{equation}
for $J = \{1,2,\ldots,k\}$ from Definition \ref{def:contrast-matrix}. 

To show $(\ref{eqn:cond-for-orthogonal-y}) \Rightarrow
(\ref{eqn:cond-for-y_J-orthogonal})$, we use an induction for $t$. For
the case of $t=1$, we have
\[
 y_{\{j\}}(1) =  y_{\{j\}}(2) = \cdots =  y_{\{j\}}(r_j),\ \ j = 1,\ldots,n 
\]
from $C_1\By = \Bzero_{v_1}$. From $s = \displaystyle\sum_{i_j =
1}^{r_j}y_{\{j\}}(i_j)$, we have
\[
 y_{\{j\}}(i_j) = \frac{s}{r_j} = \frac{s}{m_{\{j\}}}
\]  
for $j = 1,\ldots,n$. Next consider the case of $t$ under the assumption
that the theorem holds for the case of $t-1$. 
We write $J = \{1,2,\ldots,t\}$ 
for a $J$ with $\#J = t$ without loss of generality. Our purpose is to
show that
\[
 \By_J(\Bi_J) = \frac{s}{m_J},\ \ \forall \Bi_J \in {\cal I}_J
\]
for $J = \{1,2,\ldots,t\}$. Similarly to the relation
(\ref{eqn:y_J-edge}), for $\Bi_J = (i_1,\ldots,i_{t-1},i_t) \in
\left(
\displaystyle\prod_{j = 1}^{t-1}[r_j - 1]\right) \times [r_t]$ we have
\[
\frac{s}{m_{J \setminus\{t\}}} = 
 y_{J \setminus \{t\}}(i_1,\ldots,i_{t-1}) = \sum_{i_t =
 1}^{r_t}y_J(i_1,\ldots,i_{t-1},i_t) 
\]
and therefore
\[
 y_J(i_1,\ldots,i_{t-1},i_t) = \frac{1}{r_t}\cdot \frac{s}{m_{J\setminus\{t\}}}
 = \frac{s}{m_J}
\]
holds for $i_t = 1,\ldots,r_t$. To show the relation for other
$\Bi_J$'s, suppose $p$ elements of $\{i_1,\ldots,i_{t-1}\}$ equal to
$\{r_1,\ldots,r_{t-1}\}$, i.e.,
\[
 p = \#\{i_j\ :\ i_j = r_j,\ j = 1,\ldots,t-1\}.
\]
Again we use an induction for $p$ here. For the case of $p = 1$, suppose
$i_1 = r_1$ without loss of generality. We have
\[
 y_J(r_1,i_2,\ldots,i_t) = \displaystyle
  y_{J\setminus\{1\}}(i_2,\ldots,i_t) -
 \sum_{i_1 = 1}^{r_1 - 1}y_J(i_1,i_2,\ldots,i_t).
\]
From the assumption of the induction for $t$, we have
$y_{J\setminus\{1\}}(i_2,\ldots,i_t) = \displaystyle \frac{s}{m_{J
\setminus\{1\}}}$. Also from the assumption of the induction for $p$, we
have
$y_J(i_1,i_2,\ldots,i_t) = \displaystyle\frac{s}{m_J}$ for $i_1 =
1,\ldots,r_1 - 1$. Therefore we have
\[
 y_J(r_1,i_2,\ldots,i_t) = \displaystyle \frac{s}{m_{J\setminus\{1\}}} -
 (r_1-1)\cdot \frac{s}{m_J} = \frac{s}{m_J}.
\] 
The case of $p$ under the assumption that the relation holds for the
case of $p-1$ can be shown similarly. 
\hspace*{\fill}$\Box$

\bigskip

From Theorem \ref{prop:orthogonal} and the relation $\Btheta =
X^{-1}\By$, the 
constraints to be added to the relation (\ref{eq:prop-theta-mu-eq}) for
the orthogonal designs of strength $t$\ ($t \leq n$) becomes
\[
\Bone_m^TX\Btheta = s,\ C_{\ell}X\Btheta = \Bzero_{v_{\ell}},\ \ \ell =
1,\ldots,t.
\]
This is a generalization of relation for two-level case such as
(\ref{eqn:example-constraints-2-level}).

\begin{example}[Continuation of Example \ref{example:2x2x3-3}]
Consider $\{-1,1\}\times \{-1,1\}\times \{-1,0,1\}$ designs. In addition
 to the polynomial equations derived in Example \ref{example:2x2x3-3},
the coefficients of the indicator functions of designs with size $s$
 satisfy the relation
\[
 12\theta_{000} + 8\theta_{002} = s.
\]
The constraints for the equireplicated designs, i.e., 
orthogonal designs of strength $1$, are
\[\begin{array}{l}
-12\theta_{100} -8\theta_{102} = 0,\\
-12\theta_{010} -8\theta_{012} = 0,\\
-4\theta_{001} + 4\theta_{002} = 0,\\
-8\theta_{001} = 0. 
\end{array} 
\]
Therefore for a given $s$ that is a common multiple of $\{2,2,3\}$,
 i.e., only $s = 6$ is the compatible size of the fractional factorial
 designs in this
 case,  we can enumerate all the equireplicated
 designs as the solutions of a system of these polynomial equations.
\hspace*{\fill}$\Box$
\end{example}

\bigskip

Considering Theorem \ref{prop:orthogonal} for the case of $t = n$, we
have the following.

\begin{corollary}
The contrast matrix $C$ is a non-singular $m\times m$ matrix.
\end{corollary}

Now we give another representation of the indicator function reflecting the
orthogonality. 
For the indicator function
(\ref{eqn:poly-inter-polatory}), consider a non-singular linear
transformation $\Btheta \mapsto \Bmu = CX\Btheta$. New variables $\Bz$
is also defined by $\Bz =
((CX)^{-1})^T\Bx$, 
where 
$\Bx = (\Bx^{\Ba})_{\Ba
\in L}$ is an $m\times 1$ column vector of variables, and 
$\Bz =
\{z_{J(\tilde{\Bi})}\ :\ J \subset [n], \tilde{\Bi} \in \prod_{j \in
J}[r_j - 1]\}$ is also an $m\times 1$ column vector of variables.
Then we have a representation of the indicator function 
for $\Bz$, 
\begin{equation}
 f(\Bz) = \sum_{J \subset [n], \tilde{\Bi} \in \prod_{j \in
J}[r_j - 1]} \mu_{J(\tilde{\Bi})}z_{J(\tilde{\Bi})}.
\label{eqn:contrast-expression}
\end{equation}
We call (\ref{eqn:contrast-expression})
 {\it the contrast representation} of the indicator function.

From the contrast representation, we see the size and the orthogonality of
the designs directly, which is 
the advantage of the contrast representation. For example, the constant term
$\mu_{\emptyset}$ is the size of the design, and 
\[
 \mu_{J(\tilde{\Bi})} = 0\ \ \mbox{for}\ \#J = 1, 
\tilde{\Bi} \in \prod_{j \in
J}[r_j - 1]
\]
corresponds to equireplicated designs, and so on. 

\begin{example}\label{example:F4}
In Section \ref{sec:intro}, we see the indicator function of $3^{4-2}$
 regular fractional factorial design $F_3$ in Figure
 \ref{fig:example-F3} is (\ref{eqn:indicator-function-F3}). 
The contrast representation of $F_3$ is 
\[\begin{array}{rcl}
 f(\Bz) & = & 9 + z_{123(111)} + z_{123(112)} - 
z_{123(122)} - z_{123(212)} - z_{123(221)} - z_{124(111)}\\
& & {} - z_{124(122)} + z_{124(211)} + z_{124(212)} - z_{124(221)} -
 z_{134(111)} + z_{134(121)}\\
& & {}  + z_{134(122)} - z_{134(212)} - z_{134(221)}
- z_{234(111)}
- z_{234(122)}
+ z_{234(211)}\\
& & {} + z_{234(212)}
- z_{234(221)}
- z_{1234(1111)}
- z_{1234(2221)}.
\end{array}
\]
From this representation, we see that the size of $F_3$ is $9$, and
 $F_3$ is an orthogonal design of strength $2$. 

Another example is a
 $1/2$ fraction of $2\times 2\times 3$ design $F_4$ displayed in Figure
 \ref{fig:2x2x3-half}. 
\begin{figure*}[htbp]
\[
\begin{array}{|rrr|}
\multicolumn{3}{l}{F_4}\\
\multicolumn{1}{c}{x_1} & x_2 & \multicolumn{1}{c}{x_3}\\ \hline
-1 &-1 &-1 \\
-1 &-1 & 1 \\
-1 & 1 & 0 \\
 1 &-1 & 0 \\
 1 &-1 & 1 \\
 1 & 1 &-1 \\ \hline
 \end{array}
\]
\caption{An example of $1/2$ fraction of $\{-1,1\}^2\times\{-1,0,1\}$ design.} 
\label{fig:2x2x3-half}. 
\end{figure*}
The indicator function and the contrast representation of $F_4$ are
\begin{equation}
 f(x_1,x_2,x_3) = \frac{1}{2} - \frac{1}{2}x_1x_2 - \frac{1}{4}x_2x_3 -
 \frac{1}{4}x_1x_2x_3 - \frac{1}{4}x_2x_3^2 + \frac{3}{4}x_1x_2x_3^2
\label{eqn:indicator-function-F4} 
\end{equation}
and
\begin{equation}
 f(\Bz) = 6 + 2z_{2(1)} + z_{12(11)} - z_{23(12)} + z_{123(111)},
\label{eqn:contrast-representation-F4}
\end{equation}
respectively. From the contrast representation, we see that the size of
 $F_4$ is $6$. We also see that $x_1$ and $x_3$ are orthogonal from
\[
 \mu_{1(1)} = \mu_{3(1)} = \mu_{3(2)} = \mu_{13(11)} = \mu_{13(12)} = 0.
\] 
On the other hand, $\mu_{2(1)} \neq 0$ implies that $F_3$ is not
 equireplicated for $x_2$. 
\hspace*{\fill}$\Box$
\end{example}

In addition, 
the contrast representation does not depend on the level
coding, 
whereas the indicator function depends on the level coding. This
is another advantage of the contrast representation.

\begin{proposition}
In the contrast representation (\ref{eqn:contrast-expression}), $\Bmu$
is determined only from the contrast matrix $C$. 
It does not depend on the model matrix $X$, especially on the level-coding.
\end{proposition}

\paragraph*{Proof.}\ From  $\Btheta = X^{-1}\By$, we have $\Bmu =
 CX\Btheta = C\By$.\hspace*{\fill}$\Box$

\bigskip

In other words, the influence of the level-coding on the contrast
representation is involved in the variables $\Bz$. For the same contrast
matrix $C$ and the response $\By \in \{0,1\}^m$ given by 
 (\ref{eqn:y-subset-design}), the contrast representation of the design
 $F = \{\Bd_i\in D\ :\ \Bi \in {\cal I}'\}$ has the same coefficient vector
 $\Bmu = C\By$ regardless of the level-coding. On the other hand, the
 variable $\Bz$ depends on the level-coding and 
is defined by $\Bz = ((CX)^{-1})^{T}\Bx$.

\bigskip

Solving a system of polynomial equations for the coefficients of
the indicator function or the contrast
representation by computational algebraic softwares, 
we can obtain the complete list of fractional factorial designs with
given orthogonality in theory. 
It is true that the computational feasibility is an important issue,
which we see in Section 4. Another important point arises in classifying
the solutions to the 
equivalence classes for permutations of levels or factors.
For two-level cases, as we see in
\cite{Fontana-Pistone-Rogantin-2000}, the 
equivalence classes 
for permutations of levels and factors are simply obtained by
sign changes or permutation of indices 
for the coefficients of the indicator functions. 
To consider multi-level cases, we give the description of 
the equivalence classes as follows. Suppose $S_{{\cal I}}$ is a
group of permutations of ${\cal I}$, and $G \subset S_{{\cal I}}$ is a
group we consider, i.e., a group of permutations of levels for each factor and
permutations of factors if possible. For each $g \in
G$, let $P_g$ be an $m \times m$ permutation matrix. 
Then we have the following.

\begin{proposition}\label{prop:invariance-equivalence-class}
Let $G \subset S_{{\cal I}}$ is a group. Then  
the equivalence classes for $\Btheta$ and $\Bmu$ with respect to $G$ 
are
\[
 [\Btheta] = \{X^{-1}P_gX\Btheta\ :\ g \in G\}
\]
and
\[
 [\Bmu] = \{CP_gC^{-1}\Bmu\ :\ g \in G\},
\]
respectively.
\end{proposition}

\paragraph*{Proof.}\ Let $\tilde{\By} = P_g\By$. Let the corresponding
indicator functions be $f(\Bx) = \Btheta^{T}\Bx$ and $\tilde{f}(\Bx) =
\tilde{\Btheta}^{T}\Bx$, where $\Btheta = X^{-1}\By$ and
$\tilde{\Btheta} = X^{-1}\tilde{\By}$, respectively.
Then we have the relation
\[
 \tilde{\Btheta} = X^{-1}\tilde{\By} = X^{-1}P_g\By = X^{-1}P_gX\Btheta.
\]
Similarly, for the constant representations 
$f(\Bz) = \Bmu^{T}\Bz$ and $\tilde{f}(\Bz) = \tilde{\Bmu}^{T}\Bz$ where
$\Bmu = C\By$ and $\tilde{\Bmu} = C\tilde{\By}$, respectively, we have
\[
\tilde{\mu} = C\tilde{\By} = CP_g\By = CP_gC^{-1}\Bmu.
\] 
\hspace*{\fill}$\Box$

\bigskip

Proposition \ref{prop:invariance-equivalence-class} shows that neither the
indicator function nor the contrast representation has the
invariance property for multi-level cases. We will see it in the
computations in Section
\ref{sec:computation}.

\section{Classifications of orthogonal $2^3\times 3$ and $2^4\times 3$ designs}
\label{sec:computation}

In this section, we consider $2^3 \times 3$ and $2^4\times 3$
designs. Using a computational algebraic software, we solve systems of the
polynomial equations and derive a classification of designs with given
characteristic. All the computations are done by Macaulay2
(\cite{Macaulay2}) installed in a virtual machine (vmware) 
on a laptop with $2.80$ GHz CPU and $8$ GB memory. The memory allocated to
the virtual machine is $512$ MB.

\subsection{Full enumeration of the orthogonal fractions of the $2^3\times 3$
  designs of strength $2$}
First we consider the orthogonal fractions of the $2^3\times 3$ designs
of strength $2$. Corresponding system of algebraic equations includes a
set of $m = 2^3 \times 3 = 24$ general relations, $1$ relation for the size,
$5$ relations for the balance for each factor and $9$ relations
for the orthogonality of strength $2$, for $25$ variables. Note that
there are $m + 1$ variables, where $+1$ corresponds to the variable for
the size $s$. To obtain a compatible size $s$, first we calculate the
Gr\"obner basis of the ideal $I$ generated by the $39$ polynomials
corresponding to the above $39$ relations 
for the elimination ordering where the variable $s$ is the lowest. 
For the level-coding $\{-1,1\}^3\times \{-1,0,1\}$, 
the Gr\"obner basis is calculated within $0.1$ seconds, and the
elimination ideal is
\[
 I \cap \mathbb{Q}[s] = \left<\ s^3 - 36s^2 + 288s\ \right>
= \left<\ s(s - 12)(s - 24)\ \right>,
\] 
i.e., only the size $s = 12$ is compatible. This result is also obvious
because the size of the orthogonal designs must be the multiple of
$2\times 2$ and $2\times 3$. 
Note that the Gr\"onber basis calculations heavily depend on the
level-coding. To
see this, the author also try the same computation under the level-coding
$\{0,1\}^3\times \{0,1,2\}$, and find that the computation does not
finish in one week.

Now we fix $s = 12$ and calculate all the solutions. We find that
there are $44$ solutions, classified into $3$ equivalence classes as
follows.
\begin{itemize}
\item Type (a): $2$ relations. 
The indicator function and the contrast
      representation of the representative fraction displayed in Figure 
\ref{fig:res-2x2x2x3}(a) are 
\[
 \frac{1}{2} + \frac{1}{2}x_1x_2x_3
\] 
and 
\[
 12 - 3z_{123(111)}, 
\]
respectively. This is a class of the regular fractional factorial
      designs with the defining relation $x_1x_2x_3 = 1$.
\item Type (b): $6$ relations. The indicator function and the contrast
      representation of the representative fraction displayed in Figure 
\ref{fig:res-2x2x2x3}(b) are 
\[
 \frac{1}{2} + \frac{1}{2}x_1x_2x_3 - \frac{1}{2}x_1x_2x_3x_4 -
      \frac{1}{2}x_1x_2x_3x_4^2
\] 
and 
\[
 12 - z_{123(111)} - z_{1234(1112)}, 
\]
respectively, each with $4$ relations. The indicator function and
      the contrast representation of another fraction in the same
      equivalence class are
\[
 \frac{1}{2} - \frac{1}{2}x_1x_2x_3 + x_1x_2x_3x_4^2
\] 
and 
\[
 12 - z_{123(111)} - z_{1234(1111)}, 
\]
respectively, each with $2$ relations.
\item Type (c): $36$ relations. The indicator function and the contrast
      representation of the representative fraction displayed in Figure 
\ref{fig:res-2x2x2x3}(c) are 
\[
 \frac{1}{2} + \frac{1}{2}x_1x_2x_3 - \frac{1}{2}x_1x_3x_4 -
      \frac{1}{2}x_1x_2x_3x_4^2 
\] 
and 
\[
 12 - z_{123(111)} + z_{134(111)} + 2z_{134(112)} + z_{1234(1111)} +
      z_{1234(1112)}, 
\]
respectively, each with $12$ relations. The indicator function and
      the contrast representation of another fraction in the same
      equivalence class are
\[
 \frac{1}{2} - \frac{1}{2}x_1x_2 - \frac{1}{4}x_1x_2x_4 +
      \frac{1}{4}x_1x_2x_3x_4 + \frac{3}{4}x_1x_2x_4^2 +
      \frac{1}{4}x_1x_2x_3x_4^2 
\] 
and 
\[
 12 - z_{123(111)} + 2z_{124(111)} + z_{124(112)} + z_{1234(1111)} +
      z_{1234(1112)}, 
\]
respectively, each with $24$ relations.
\end{itemize} 
\begin{figure*}[htbp]
\[
\begin{array}{|rrrr|}
\multicolumn{4}{l}{{\rm Type (a)}}\\
\multicolumn{1}{c}{x_1} & x_2 & x_3 & \multicolumn{1}{c}{x_4}\\ \hline
 1 & 1 & 1 &-1\\
 1 & 1 & 1 & 0\\
 1 & 1 & 1 & 1\\
 1 &-1 &-1 &-1\\
 1 &-1 &-1 & 0\\
 1 &-1 &-1 & 1\\
-1 & 1 &-1 &-1\\
-1 & 1 &-1 & 0\\
-1 & 1 &-1 & 1\\
-1 &-1 & 1 &-1\\
-1 &-1 & 1 & 0\\
-1 &-1 & 1 & 1\\ \hline
 \end{array}
\hspace*{10mm}
\begin{array}{|rrrr|}
\multicolumn{4}{l}{{\rm Type (b)}}\\
\multicolumn{1}{c}{x_1} & x_2 & x_3 & \multicolumn{1}{c}{x_4}\\ \hline
 1 & 1 & 1 &-1\\
 1 & 1 & 1 & 0\\
 1 & 1 &-1 & 1\\
 1 &-1 & 1 & 1\\
 1 &-1 &-1 & 0\\
 1 &-1 &-1 &-1\\
-1 & 1 & 1 & 1\\
-1 & 1 &-1 & 0\\
-1 & 1 &-1 &-1\\
-1 &-1 & 1 &-1\\
-1 &-1 & 1 & 0\\
-1 &-1 &-1 & 1\\ \hline
 \end{array}
\hspace*{10mm}
\begin{array}{|rrrr|}
\multicolumn{4}{l}{{\rm Type (c)}}\\
\multicolumn{1}{c}{x_1} & x_2 & x_3 & \multicolumn{1}{c}{x_4}\\ \hline
 1 & 1 & 1 &-1\\
 1 & 1 & 1 & 0\\
 1 & 1 &-1 & 1\\
 1 &-1 & 1 &-1\\
 1 &-1 &-1 & 0\\
 1 &-1 &-1 & 1\\
-1 & 1 & 1 & 1\\
-1 & 1 &-1 & 0\\
-1 & 1 &-1 &-1\\
-1 &-1 & 1 & 1\\
-1 &-1 & 1 & 0\\
-1 &-1 &-1 &-1\\ \hline
 \end{array}
\]
\caption{Orthogonal fractions of the $2^3\times 3$ designs of strength $2$}
\label{fig:res-2x2x2x3}
\end{figure*}

In the above list, Type (a) is the class of the regular fractions,
whereas Type (b) and Type (c) are classes of the non-regular
fractions. 
Note that Type (b) and Type (c) differ only in the last columns (the
levels of $x_4$) in Figure \ref{fig:res-2x2x2x3}. For each row where the
levels of $(x_1,x_2,x_3)$ is unique, there are $4$ such rows, the levels
of $x_4$ are fixed ($x_4 = 1$) in Type (b), whereas the levels of $x_4$
are $1$ or $-1$ in Type (c).
Type (b) and (c) can also be characterized considering 
the designs
obtained from a traditional $OA(12, 3^1 2^4)$ orthogonal array as
follows. The orthogonal array $OA(12, 3^1 2^4)$ in Appendix 8C of 
\cite{Wu-Hamada-2009} is displayed in Figure \ref{OA12}.
\begin{figure*}
\[
 \begin{array}{ccccc}\hline
1 & 2 & 3 & 4 & 5\\ \hline
0 & 0 & 0 & 0 & 0\\
0 & 0 & 1 & 0 & 1\\
0 & 1 & 0 & 1 & 1\\
0 & 1 & 1 & 1 & 0\\
1 & 0 & 0 & 1 & 1\\
1 & 0 & 1 & 1 & 0\\
1 & 1 & 0 & 0 & 1\\
1 & 1 & 1 & 0 & 0\\
2 & 0 & 0 & 1 & 0\\
2 & 0 & 1 & 0 & 1\\
2 & 1 & 0 & 0 & 0\\
2 & 1 & 1 & 1 & 1\\ \hline
 \end{array}
\]
\caption{Orthogonal array $OA(12, 3^1 2^4)$  (Appendix 8C of
 \cite{Wu-Hamada-2009}).}
\label{OA12}
\end{figure*}
From $OA(12, 3^12^4)$, we can obtain $1/2$ fractions of $2^3\times 3$
designs by selecting $3$ columns from the columns $\{2,3,4,5\}$. 
We see that all the designs
constructed in this way are included in the equivalence class of Type
(c). Therefore Type (c) is regarded as the class of $OA(12, 3^1 2^4)$ designs.

\subsection{Full enumeration of the orthogonal fractions of the $2^4\times 3$
  designs of strength $3$}
Next we consider the fractions of the $2^4\times 3$ designs. 
For this case, enumeration of the orthogonal fractions of strength $2$
may be difficult to compute for standard PC. In fact, the Gr\"obner
basis of the elimination ideal for the compatible size does
not obtained after $1$ week calculation under the level-coding
$\{-1,1\}^4\times \{-1,0,1\}$. 
Therefore we enumerate the
orthogonal fractions of strength $3$ instead. Note that, for fixed size $s$,
there are $m = 48$ variables with constraints $1 + 6 + 14 = 21$
relations for strength $2$, and with constraints $1 + 6 + 14 + 16 = 37$
relations for strength $3$. Therefore, by eliminating variables, the
number of variables reduces 
$11$ for strength $3$, whereas to $27$ for strength $2$. The compatible
size must be $s = 24$ for strength $3$, that is only the multiple of $2\times
2\times 2$ and $2\times 2\times 3$ less than $m = 48$. This fact is also
checked by the Gr\"obner basis calculation. After calculation within
$0.1$ seconds, we see that the elimination ideal is
\[
 I\cap \mathbb{Q}[s] = \left<\ s^3 - 72s^2 + 1152s\ \right>
= \left<\ s(s - 24)(s - 48)\ \right>.
\]
Therefore we fix $s = 24$ and calculate all the solutions. We find there
are $56$
solutions, classified into $3$ equivalence classes as follows.
\begin{itemize}
\item Type (a): 2 relations. The indicator function and the contrast
      representation of the representative fraction displayed in Figure
      \ref{fig:res-2x2x2x2x3}(a) are
\[
 \frac{1}{2} + \frac{1}{2}x_1x_2x_3x_4
\]
and
\[
 24 + 3z_{1234(1111)},
\]
respectively. This is a class of the regular fractional factorial
      designs with the defining relation $x_1x_2x_3x_4 = 1$.

\item Type (b): 6 relations. 
The indicator function and the contrast
      representation of the representative fraction displayed in Figure
      \ref{fig:res-2x2x2x2x3}(b) are
\[
 \frac{1}{2} - \frac{1}{2}x_1x_2x_3x_4 - \frac{1}{2}x_1x_2x_3x_4x_5 +
      \frac{1}{2}x_1x_2x_3x_4x_5^2
\]
and
\[
 24 - z_{1234(1111)} + z_{12345(11111)} + z_{12345(11112)},
\]
respectively, each with 4 relations. 
The indicator function and the
      contrast representation of another fraction in the same
      equivalence class are
\[
 \frac{1}{2} + \frac{1}{2}x_1x_2x_3x_4 - x_1x_2x_3x_4x_5^2
\]
and
\[
 24 - z_{1234(1111)} - z_{12345(11111)},
\]
respectively, each with 2 relations.

\item Type (c): 48 relations. The indicator function and the contrast
      representation of the representative fraction displayed in Figure
      \ref{fig:res-2x2x2x2x3}(c) are
\[
 \frac{1}{2} -  \frac{1}{2}x_1x_2x_3x_4 -  \frac{1}{2}x_1x_2x_4x_5 +
      \frac{1}{2}x_1x_2x_3x_4x_5^2
\]
and
\[
 24 - z_{1234(1111)} - z_{1245(1111)} - 2z_{1245(1112)} - z_{12345(11112)},
\]
respectively, each with 16 relations. 
The indicator function and the
      contrast representation of another fraction in the same
      equivalence class are
\[
 \frac{1}{2} + \frac{1}{2}x_1x_2x_3 + \frac{1}{4}x_1x_2x_3x_5 +
      \frac{1}{4}x_1x_2x_3x_4x_5 - \frac{3}{4}x_1x_2x_3x_5^2 +
      \frac{1}{4}x_1x_2x_3x_4x_5^2
\]
and
\[
 24 + z_{1234(1111)} + 2z_{1235(1111)} + z_{1235(1112)} + z_{12345(11111)},
\]
respectively, each with 32 relations.
\end{itemize}

\begin{figure*}[htbp]
\[
\begin{array}{|rrrrr|}
\multicolumn{5}{l}{{\rm Type (a)}}\\
\multicolumn{1}{c}{x_1} & x_2 & x_3 & x_4 & \multicolumn{1}{c}{x_5}\\ \hline
-1 & -1 & -1 & -1 & -1\\
-1 & -1 & -1 & -1 &  0\\
-1 & -1 & -1 & -1 &  1\\
-1 & -1 &  1 &  1 & -1\\
-1 & -1 &  1 &  1 &  0\\
-1 & -1 &  1 &  1 &  1\\
-1 &  1 & -1 &  1 & -1\\
-1 &  1 & -1 &  1 &  0\\
-1 &  1 & -1 &  1 &  1\\
-1 &  1 &  1 & -1 & -1\\
-1 &  1 &  1 & -1 &  0\\
-1 &  1 &  1 & -1 &  1\\
 1 & -1 & -1 &  1 & -1\\
 1 & -1 & -1 &  1 &  0\\
 1 & -1 & -1 &  1 &  1\\
 1 & -1 &  1 & -1 & -1\\
 1 & -1 &  1 & -1 &  0\\
 1 & -1 &  1 & -1 &  1\\
 1 &  1 & -1 & -1 & -1\\
 1 &  1 & -1 & -1 &  0\\
 1 &  1 & -1 & -1 &  1\\
 1 &  1 &  1 &  1 & -1\\
 1 &  1 &  1 &  1 &  0\\
 1 &  1 &  1 &  1 &  1\\ \hline
 \end{array}
\hspace*{10mm}
\begin{array}{|rrrrr|}
\multicolumn{5}{l}{{\rm Type (b)}}\\
\multicolumn{1}{c}{x_1} & x_2 & x_3 & x_4 & \multicolumn{1}{c}{x_5}\\ \hline
-1 & -1 & -1 & -1 & -1\\
-1 & -1 & -1 &  1 &  0\\
-1 & -1 & -1 &  1 &  1\\
-1 & -1 &  1 & -1 &  1\\
-1 & -1 &  1 & -1 &  0\\
-1 & -1 &  1 &  1 & -1\\
-1 &  1 & -1 & -1 &  1\\
-1 &  1 & -1 & -1 &  0\\
-1 &  1 & -1 &  1 & -1\\
-1 &  1 &  1 & -1 & -1\\
-1 &  1 &  1 &  1 &  0\\
-1 &  1 &  1 &  1 &  1\\
 1 & -1 & -1 & -1 &  1\\
 1 & -1 & -1 & -1 &  0\\
 1 & -1 & -1 &  1 & -1\\
 1 & -1 &  1 & -1 & -1\\
 1 & -1 &  1 &  1 &  0\\
 1 & -1 &  1 &  1 &  1\\
 1 &  1 & -1 & -1 & -1\\
 1 &  1 & -1 &  1 &  0\\
 1 &  1 & -1 &  1 &  1\\
 1 &  1 &  1 & -1 &  1\\
 1 &  1 &  1 & -1 &  0\\
 1 &  1 &  1 &  1 & -1\\ \hline
 \end{array}
\hspace*{10mm}
\begin{array}{|rrrrr|}
\multicolumn{5}{l}{{\rm Type (c)}}\\
\multicolumn{1}{c}{x_1} & x_2 & x_3 & x_4 & \multicolumn{1}{c}{x_5}\\ \hline
-1 & -1 & -1 & -1 &  1\\
-1 & -1 & -1 &  1 &  0\\
-1 & -1 & -1 &  1 & -1\\
-1 & -1 &  1 & -1 &  1\\
-1 & -1 &  1 & -1 &  0\\
-1 & -1 &  1 &  1 & -1\\
-1 &  1 & -1 & -1 & -1\\
-1 &  1 & -1 & -1 &  0\\
-1 &  1 & -1 &  1 &  1\\
-1 &  1 &  1 & -1 & -1\\
-1 &  1 &  1 &  1 &  0\\
-1 &  1 &  1 &  1 &  1\\
 1 & -1 & -1 & -1 & -1\\
 1 & -1 & -1 & -1 &  0\\
 1 & -1 & -1 &  1 &  1\\
 1 & -1 &  1 & -1 & -1\\
 1 & -1 &  1 &  1 &  0\\
 1 & -1 &  1 &  1 &  1\\
 1 &  1 & -1 & -1 &  1\\
 1 &  1 & -1 &  1 &  0\\
 1 &  1 & -1 &  1 & -1\\
 1 &  1 &  1 & -1 &  1\\
 1 &  1 &  1 & -1 &  0\\
 1 &  1 &  1 &  1 & -1\\ \hline
 \end{array}
\]
\caption{Orthogonal fractions of the $2^4\times 3$ designs of strength $3$}
\label{fig:res-2x2x2x2x3}
\end{figure*}

An interpretation of this list is similar to the $2^3\times 3$ case. In
Figure \ref{fig:res-2x2x2x2x3}, Type (b) and Type (c) differ only in the
last column (the levels of $x_5$). For each row where the levels of
$(x_1,x_2,x_3,x_4)$ is unique, there are $8$ such rows, the levels of
$x_5$ are fixed ($x_5 = -1$) in Type (b), whereas the levels of $x_5$
are $1$ or $-1$ in Type (c).

\section{Discussion}
In this paper, we give how to construct a system of polynomial equations
for the coefficients of the indicator functions of multi-level
fractional factorial designs with given orthogonality. We also define
the contrast representation of the indicator function, which reflects
the orthogonality of the design directly. 
The contrast representation has an advantage that it does not depends on
the level-coding.
Using these results, we show
the classifications of the orthogonal fractions of the $2^3\times 3$
designs with strength $2$ 
and $2^4\times 3$ designs with strength $3$.

In theory, we can obtain classifications of fractional factorial
designs for any size by our method. However, the computational
feasibility depends on the size of problems. For the class of 
$2^m \times 3$ designs, we see that the $2^4\times 3$ problem for
orthogonality with strength $3$ is easy to calculate. However, a 
$2^5\times 3$ problem of strength $3$ orthogonality seems very difficult
to compute. As for the class of $2^m \times 3^2$ designs, we
find that the
$2^3 \times 3^2$ problem of strength $2$ orthogonality is hard to
compute, i.e., the Gr\"obner basis calculation for the elimination ideal
does not finish in $1$ week. In addition, it is obvious that there is
no orthogonal fractions $2^3 \times 3^2$ with strength $3$ because the
size must be a multiple of $2\times 2\times 2$ and $2\times 3\times 3$. 
As a consequence, we only have limited 
computational results in this paper.

In particular, the merit of the contrast representation must be
investigated from the computational aspects. It seems that 
a system of polynomial equations for $\Bmu$ is easy to solve than that for
$\Btheta$. Note that 
the polynomial relations for $\Bmu$ are obtained by substituting
$\Btheta = (CX)^{-1}\Bmu$ to the polynomial relations for
$\Btheta$. Therefore, translating the relations for $\Btheta$ to the
relations for $\Bmu$ corresponds to the matrix operations of inverse in
advance. Unfortunately, for the problems considered in this paper, the
systems of
the polynomial equations for $\Btheta$ and $\Bmu$ are both quite easy or
quite difficult, and the effectiveness of the contrast representation
from the computational aspect is not shown. 
Therefore the quantitative evaluation of the effect of this
transformation is one of the open problems. 
It is also an open problem to compare our method to the brute-force search.

\bibliographystyle{plain}

\begin{thebibliography}{99}
\bibitem{Aoki-Takemura-2009}
S. Aoki and A. Takemura. (2009).
Some characterizations of affinely full-dimensional factorial
	designs. {\it Journal of Statistical
        Planning and Inference}, {\bf 139}, 3525--3532.
\bibitem{Box-Hunter-1961}
G. E. P. Box and J. S. Hunter. (1961).
The $2^{k-p}$ fractional factorial design. {\it Technometrics}, {\bf 3},
	311--351, 449--458.
\bibitem{Carlini-Pistone-2009}
E. Carlini and G. Pistone. (2009).
Hilbert bases for orthogonal arrays. {\it Journal of Statistical Theory
	and Practice}, {\bf 1}, 299-309.
\bibitem{Cheng-Yw-2004}
S. W. Cheng and K. Q. Ye. (2004).
Geometric isomorphism and minimum aberration for factorial designs with
	quantitative factors. {\it The Annals of Statistics}, {\bf 32},
	2168--2185.
\bibitem{Cox-Little-OShea-1992}
D. Cox, J. Little, and D. O'Shea. (2007).
{\it Ideals, varieties, and algorithms, An introduction to computational
	     algebraic geometry and commutative algebra}, Third ed.,
 Springer-Verlag.
\bibitem{Fontana-Pistone-Rogantin-2000}
R. Fontana, G. Pistone and M. P. Rogantin. (2000). Classification of
        two-level factorial fractions. {\it Journal of Statistical
        Planning and Inference}, {\bf 87}, 149--172.
\bibitem{Lauritzen}
S. L. Lauritzen. {\it Graphical models}. (1996). Oxford Statistical Science
	Series, Oxford Science Publications. 
\bibitem{Macaulay2}
D. R. Grayson and M. E. Stillman. 
Macaulay2, a software system for research in algebraic geometry.
Available at {\tt http://www.math.uiuc.edu/Macaulay2/}.
\bibitem{dojo-en}
T. Hibi (ed.) (2013).
{\it Gr\"obner bases, Statistics and software systems}. Springer.
\bibitem{Pistone-Riccomagno-Wynn-2001}
G. Pistone, E. Riccomagno and H. P. Wynn. (2001).
{\it Algebraic Statistics: Computational Commutative Algebra in Statistics}. 
Chapman \& Hall, London.
\bibitem{Pistone-Rogantin-2008b}
G. Pistone and M. P. Rogantin. (2008).
Indicator function and complex coding for mixed fractional factorial
       designs. {\it Journal of Statistical Planning and Inference},
       {\bf 138}, 787--802.
\bibitem{Pistone-Wynn-1996}
G. Pistone and H. P. Wynn. (1996).
Generalised confounding with Gr\"obner bases. {\it Biometrika}, {\bf
        83}, 653--666.
\bibitem{Wu-Hamada-2009}
C. F. Jeff Wu and M. S. Hamada. (2009).
{\em Experiments: Planning, analysis, and parameter design
  optimization}. 2nd ed. 
Wiley Series in Probability and Statistics: Texts and References
  Section. John Wiley \& Sons Inc., New York.
A Wiley-Interscience Publication.
\bibitem{Ye-2003}
K. Q. Ye. (2003). Indicator function and its application in two-level
       factorial designs. {\it Annals of Statistics}, {\bf 31},
       984--994.

\end{thebibliography}

\end{document}